# Comment

## Peter L. Bartlett, Michael I. Jordan and Jon D. McAuliffe

### INTRODUCTION

The support vector machine (SVM) has played an important role in bringing certain themes to the fore in computationally oriented statistics. However, it is important to place the SVM in context as but one member of a class of closely related algorithms for nonlinear classification. As we discuss, several of the "open problems" identified by the authors have in fact been the subject of a significant literature, a literature that may have been missed because it has been aimed not only at the SVM but at a broader family of algorithms. Keeping the broader class of algorithms in mind also helps to make clear that the SVM involves certain specific algorithmic choices, some of which have favorable consequences and others of which have unfavorable consequences—both in theory and in practice. The broader context helps to clarify the ties of the SVM to the surrounding statistical literature.

We have at least two broader contexts in mind for the SVM. The first is the family of "large-margin" classification algorithms—a class that includes boosting and logistic regression. All of these algorithms involve the minimization of a convex contrast or loss function that upper bounds the 0–1 loss function. The SVM makes a specific choice of convex loss function—the so-called hinge loss. Hinge loss has some potentially desirable properties (e.g., sparseness) and some potentially undesirable properties (e.g., lack of calibration to posterior probabilities). As we discuss, much of the theoretical analysis of the SVM is best carried out by focusing on convexity and abstracting away from the details of specific loss functions.

Second, as the authors note, the SVM is an instance of the broader family of statistical procedures based on reproducing kernel Hilbert spaces (RKHSs). The authors' emphasis is on the use of RKHS methods to provide basis expansions for discriminant functions and regression functions. RKHS ideas have, however, been carried significantly further in recent years, enlivening areas of computationally oriented statistics beyond classification and regression. We wish to convey some of the reasons for this broader interest in RKHS-based approaches.

There are both computational and statistical motivations for focusing on methods based on convexity and reproducing kernel Hilbert spaces. In the remainder of this discussion we attempt to disentangle some of these motivations, but we wish to emphasize at the outset that it is precisely because these methods bring computational and statistical considerations together that they are so interesting.

### CONVEXITY

The SVM is one example of a general strategy for solving the binary classification problem via a "convex surrogate loss function." To develop this perspective, let us define binary classification as the problem of choosing a discriminant function $f : \mathcal{X} \to \mathbb{R}$ that minimizes misclassification risk,

$$R(f) = P(Y \neq \operatorname{sgn}(f(X))) = \mathbf{E}\ell(Yf(X)),$$

where $X \in \mathcal{X}$ is the covariate, $Y \in \{\pm 1\}$ is the binary response, and $\ell(\alpha) = 1$ for $\alpha \leq 0$ and $= 0$ otherwise. The family of large-margin classification algorithms attacks this problem indirectly by minimizing a quantity known as the $\phi$-risk,

$$R_\phi(f) = \mathbf{E}\phi(Yf(X)),$$

where $\phi : \mathbb{R} \to \mathbb{R}$ is a surrogate for the loss function $\ell$, and $Yf(X)$ is called the *margin* of $f$ on the observation $(X, Y)$. The margin indicates not only whether the observation is correctly classified by $f$, but how close $f$ comes to choosing the opposite label. The surrogate loss function $\phi$ is chosen so that large margins correspond to small losses.

*Peter L. Bartlett and Michael I. Jordan are Professors, Computer Science Division and Department of Statistics, University of California, Berkeley, California 94720, USA e-mail: bartlett@stat.berkeley.edu; jordan@stat.berkeley.edu. Jon D. McAuliffe is Assistant Professor, Department of Statistics, The Wharton School, University of Pennsylvania, Philadelphia, Pennsylvania 19104-6340, USA e-mail: mcjon@wharton.upenn.edu.*

Given a data set $(X_1, Y_1), \ldots, (X_n, Y_n)$, we can form the *empirical $\phi$-risk*

$$\hat{R}_\phi(f) = \frac{1}{n} \sum_{i=1}^n \phi(Y_i f(X_i))$$

and attempt to minimize this quantity with respect to the discriminant function $f$. When $\phi$ is chosen as a convex function and $f$ is constrained to lie in a convex family of prediction rules, the minimization becomes a convex optimization problem. Contrast this with minimization of empirical 0–1 (misclassification) risk,

$$\min_{f \in \mathcal{F}} \hat{R}(f) = \frac{1}{n} \sum_{i=1}^n \ell(Y_i f(X_i)),$$

a problem whose exact or even approximate solution is known to be intractable for most nontrivial function classes $\mathcal{F}$ (e.g., Arora, Babai, Stern and Sweedyk, 1997).

In the case of SVMs, the convex surrogate is the hinge loss

$$\phi(\alpha) = \begin{cases} 1 - \alpha, & \text{if } \alpha \leq 1, \\ 0, & \text{otherwise}, \end{cases}$$

and the function $f$ is chosen from the RKHS $\mathcal{H}$ defined by the kernel. However, the hinge loss is not the only convex surrogate worth considering. Using the binomial deviance function as a convex surrogate and optimizing over linear functions on $R^p$ yields logistic regression. Just as with the SVM, a nonlinear version of logistic regression can be defined by making use of reproducing kernels (Zhu and Hastie, 2005). AdaBoost (Schapire and Singer, 1999) can also be interpreted as a large-margin method, with $\phi(\alpha) = e^{-\alpha}$; similar greedy ensemble methods correspond to other choices of $\phi$.

The benefits of empirical convex risk minimization are not just computational. Searching for a prediction rule which achieves a large margin on many training examples, rather than just correctly classifying them, is an implicit form of regularization. For example, certain function classes of infinite Vapnik–Chervonenkis (VC) dimension, where empirical 0–1 risk minimization does not yield good classifiers, can be used effectively in the large-margin framework (Bartlett, 1998; Schapire, Freund, Bartlett and Lee, 1998).

Taking the margin-based viewpoint highlights the important role convexity plays in the success of SVMs. On the other hand, the authors' emphasis on the need to find a differentiable or smooth formulation seems misplaced. In the differentiable case, the key property of a convex objective function $f$ is that, for any two points $x, y$ in the domain,

$$(1) \qquad f(y) \geq f(x) + \left(\frac{\partial f(x)}{\partial x}\right)^\top (y - x).$$

Thus, local behavior of $f$ (its gradient at $x$) determines a global lower bound on $f$. The existence of this bound makes possible efficient algorithms for convex optimization. However, property (1) holds in a slightly generalized form even for nondifferentiable convex functions. A *subgradient* of a convex $f$ at $x$ is a vector $g$ such that

$$f(y) \geq f(x) + g^\top (y - x) \quad \forall y.$$

The *subdifferential* $\partial f(x)$ of $f$ at $x$ is the set of $f$'s subgradients at $x$. The subdifferential is the natural analog of the gradient for nonsmooth objectives: any point in $\partial f(x)$ provides the equivalent of property (1); $0 \in \partial f(x)$ if and only if $x$ is a global minimizer of $f$; and the subdifferential contains only the gradient at points of differentiability. Moreover, every convex function has nonempty subdifferentials throughout the interior of its domain. There has been a great deal of successful research on efficient algorithms for nonsmooth convex optimization using "bundle methods" based on subdifferentials. See, for example, Hiriart-Urruty and Lemaréchal (1993), Borwein and Lewis (2000) and Boyd and Vandenberghe (2004).

**Statistical Analysis**

To address issues such as statistical consistency and finite sample behavior in the large-margin framework, we need to decompose the risk $R(f_n)$ into three components. Two of the components are the approximation error and the estimation error that are familiar from other areas of nonparametric statistics. The third component arises from the use of the convex surrogate $\phi$ in place of the 0–1 loss.

We can quantify the effect of this third component through an inequality of the form

$$\psi(R(f) - R^*) \leq R_\phi(f) - R_\phi^*,$$

where $\psi$ is a convex, nonnegative function and $f$ is an arbitrary measurable function (Bartlett, Jordan and McAuliffe, 2006). Notice that such an inequality relates the *excess risk*, $R(f) - R^*$, to the *excess $\phi$-risk*, $R_\phi(f) - R_\phi^*$. Here, the Bayes risk, $R^*$, is defined by $R^* = \inf_g R(g)$, where the infimum is over all measurable functions, and $R_\phi^*$ is the minimal $\phi$-risk, $R_\phi^* = \inf_g R_\phi(g)$. An optimal inequality of this form has been obtained for any nonnegative surrogate loss function $\phi$, where $\psi$ can be written explicitly in terms of $\phi$ (Bartlett, Jordan and McAuliffe, 2006). In the case of SVMs, where $\phi$ is the hinge loss, $\psi$ turns out to be the identity (Zhang 2004; Blanchard, Bousquet and Massart, 2006).

Thus, for SVMs we can write an optimal upper bound on the excess risk as

$$R(f_n) - R^*$$
$$\leq R_\phi(f_n) - R_\phi^*$$
$$= \left(R_\phi(f_n) - \inf_{f \in \mathcal{H}} R_\phi(f)\right) + \left(\inf_{f \in \mathcal{H}} R_\phi(f) - R_\phi^*\right).$$

This decomposition into the estimation error and the approximation error reflects the bias–variance trade-off, ubiquitous in nonparametric estimation. On the one hand, the function $f_n$ must come from a suitably simple class to ensure that the performance of $f_n$ on the finite training sample is representative of its true performance and so the estimation error $R_\phi(f_n) - \inf_{f \in \mathcal{H}} R_\phi(f)$ is not too large. On the other hand, $f_n$ must be suitably complex, so that its $\phi$-risk is not too much larger than the optimal $R_\phi^*$. One common approach to this model selection problem is the method of sieves, where the function $f_n$ is chosen as the minimizer of the empirical $\phi$-risk over classes $\mathcal{F}_n$ that grow progressively richer as the sample size $n$ increases. The other common approach is to add a regularization term, so that $f_n$ is chosen as the minimizer of

$$\hat{R}_\phi(f) + \lambda_n \Omega(f)$$

for some regularization functional $\Omega$ that penalizes complex functions and for some regularization coefficients $\lambda_n$. In the case of SVMs, the regularization functional is the squared norm in the reproducing kernel Hilbert space $\mathcal{H}$. The regularization approach and the method of sieves are closely related. In particular, since $0 \in \mathcal{H}$ and $\hat{R}_\phi(0) = 1$ for the hinge loss, we know that

$$\|f_n\|_\mathcal{H}^2 \leq \frac{1}{\lambda_n}(\hat{R}_\phi(f_n) + \lambda_n \|f_n\|_\mathcal{H}^2)$$
$$\leq \frac{1}{\lambda_n}(\hat{R}_\phi(0) + \lambda_n \|0\|_\mathcal{H}^2) = \frac{1}{\lambda_n}.$$

Thus, the SVM chooses a function from the sieve $\mathcal{F}_n = \{f \in \mathcal{H} : \|f\|_\mathcal{H}^2 \leq 1/\lambda_n\}$.

It is important to note that the regularization coefficient $\lambda_n$ must decrease with the sample size, so as to ensure universal consistency. Indeed, the approximation error, $\inf_{f \in \mathcal{F}_n} R_\phi(f) - R_\phi^*$, must go to zero. On the other hand, it should decrease sufficiently slowly that the estimation error, $R_\phi(f_n) - \inf_{f \in \mathcal{F}_n} R_\phi(f)$, also decreases to zero. (Consistency does not follow from uniform convergence of the empirical to expected error in a ball of fixed radius, as the appendix of the paper suggests.)

One of the most important consequences of the choice of a kernel is the way it affects this trade-off between the estimation and approximation errors. We can view the kernel, and hence the norm in the RKHS, as defining a complexity hierarchy. A good kernel is one for which a good approximation to $R_\phi^*$ can be obtained with a function that is low in the complexity hierarchy, in the sense that it has a small norm.

The paper lists several statistical issues as important open problems, notably the finite-sample performance of SVMs and the estimation of a posteriori class probabilities. In fact, bounds on the estimation error, $R_\phi(f_n) - \inf_{f \in \mathcal{F}_n} R_\phi(f)$, for finite sample size have been known for years. These bounds are expressed in terms of properties of the eigenvalues of the *Gram matrix*, the matrix whose entries are $k(x_i, x_j)$ (see, e.g., Shawe-Taylor, Bartlett, Williamson and Anthony, 1998; Bartlett and Shawe-Taylor, 1999; Williamson, Smola and Schölkopf, 1999; Mendelson, 2002; Bartlett, Bousquet and Mendelson, 2005; Blanchard, Bousquet and Massart, 2006). Moreover, these results provide an essential foundation for proofs of consistency (Steinwart, 2002, 2005; Zhang, 2004). In contrast, while the VC dimension is central to the analysis of methods that minimize the empirical 0–1 risk, it is not relevant to SVMs. Indeed, for any kernel that is sufficiently rich to allow a universal consistency result, the RKHS necessarily has infinite VC dimension. See Bartlett and Shawe-Taylor (1999) for an explanation of the role of the VC dimension and other complexity measures that are more appropriate to the analysis of the finite-sample behavior of SVMs.

Moreover, the problem of estimation of a posteriori class probabilities using SVM classification outputs is not open, in the following sense. It is an easy calculation to see that, for any $x$, the minimizer of the conditional expectation $\mathbf{E}[\phi(Y f(x))|X=x]$ is $f^*(x) = \text{sign}(2\eta(x) - 1)$, where $\eta(x) = \Pr[Y = 1|X = x]$. Furthermore, results of Steinwart (2003) establish that, under reasonable conditions on the kernel function and the rate at which the regularization coefficients go to zero, the function $f_n$ chosen by the SVM converges to $f^*$. Thus, asymptotically there is no information about the a posteriori class probabilities in the SVM classification outputs. Thus the SVM framework does not appear to provide an appropriate starting point for the estimation of a posteriori probabilities. This is a distinguishing feature of the hinge loss $\phi$: its minimization cannot correspond to fitting a probability model, since it is indifferent to distinct values of the class probability.

To elaborate on this point, note that one of the attractive features of SVM classifiers is their sparseness: The proportion of nonzero dual variables ("support vectors") is typically small. Steinwart (2004) has presented a beautiful relationship between the number of support vectors in an SVM and the Bayes risk. Assuming that the kernel is appropriately chosen and the regularization is reduced sufficiently slowly as the sample size increases, the asymptotic proportion of support vectors is equal to twice the Bayes risk. On the other hand, it is known that if we replace the hinge loss $\phi$ with the differentiable quadratic loss, sparseness disappears, but then the a posteriori class probabilities can be estimated asymptotically. Indeed, sparseness and the ability to estimate conditional probabilities seem to be incompatible. If the hinge loss is replaced by any of a large family of loss functions, it can be shown (Bartlett and Tewari, 2004) that the proportion of support vectors approaches the expectation of a certain function of the conditional probability $\eta(X)$, and this function is maximal for those values of $\eta(X)$ for which estimation of the conditional probability is possible asymptotically.

## KERNEL METHODS

While the SVM has helped to bring RKHS ideas to new prominence, focusing on the SVM runs the risk of limiting the appreciation of the scope and potential impact of RKHS methods. In this section we augment the presentation by Moguerza and Muñoz to provide a broader context for the understanding of the RKHS aspect of the SVM approach.

The main point that we wish to make in this section is this: A RKHS provides computationally efficient machinery for evaluating and optimizing a variety of statistical functionals of interest. The empirical loss functions of nonparametric classification and regression are a special case—one in which the focus is on the basis function expansions provided by a RKHS—but there are other roles for a RKHS.

A key idea in RKHS methodology is that an inner product can be computed by evaluating a *reproducing kernel* $k(x, x')$—a function of two arguments that obeys symmetry and positive definiteness conditions. A reproducing kernel may take the form of an analytic function (e.g., the Gaussian kernel) or may take the form of a computational procedure (e.g., the string kernel, which involves a dynamic program).

Given a set of $n$ data points $\{x_1, x_2, \ldots, x_n\}$ and given a reproducing kernel $k(x, x')$, one can form the Gram matrix. The SVM and the other kernel methods mentioned in Section 5 of Moguerza and Muñoz are all based on various operations (computation of eigenvectors, determinants, inverses, etc.) on Gram matrices. This reduction of the data to a Gram matrix is significant computationally; indeed, while the naive computational complexity of many kernel methods is $O(n^3)$, the exploitation of more sophisticated numerical linear back ends can drive this cost down to $O(nk^2)$, where $k$ is a measure of the effective rank of a Gram matrix. The effective rank is typically small.

Let us now consider a problem that at first glance seems to have little to do with kernel methods—the problem of assessing whether random variables $X_1$ and $X_2$ are independent. Independence can be reduced to correlation by considering transformations of the random variables. In particular, $X_1$ and $X_2$ are independent if and only if

$$\rho = \max_{h_1, h_2 \in \mathcal{H}} \text{Corr}(h_1(X_1), h_2(X_2)) = 0$$

for a suitably rich function space $\mathcal{H}$. Indeed, if $\mathcal{H}$ is $L_2$ and thus contains the Fourier basis, this statement reduces to a classical fact about characteristic functions. More interestingly, the result also holds for certain RKHSs. Moreover, the reproducing property of the kernel implies that function evaluation in a RKHS reduces to an inner product, $h_1(X_1) = \langle k(\cdot, X_1), h_1 \rangle$, where $k(\cdot, \cdot)$ is the reproducing kernel for $\mathcal{H}$ and $\langle \cdot, \cdot \rangle$ is the corresponding inner product. Thus correlations can be computed as

$$\begin{aligned}&\text{Corr}(h_1(X_1), h_2(X_2))\\ &\quad = \text{Corr}(\langle k(\cdot, X_1), h_1 \rangle, \langle k(\cdot, X_2), h_2 \rangle).\end{aligned}$$

Maximizing over $h_1$ and $h_2$ thus amounts to maximizing the correlation between projections of vectors in a pair of Hilbert spaces; this is nothing but canonical correlation analysis (CCA) in $\mathcal{H}$ (cf. Leurgans, Moyeed and Silverman, 1993). Moreover, using the reproducing property it is easy to show that this functional CCA computation can be reduced to a generalized eigenvector problem on a pair of Gram matrices (one Gram matrix for the observations of $X_1$ and another Gram matrix for the observations of $X_2$). Thus we can assess independence by carrying out a kernelized version of CCA.

This line of argument is due to Bach and Jordan (2002), who showed how it could be used to fit a semiparametric model known as independent component analysis. The general approach has been carried further by Gretton et al. (2005), who established relationships between RKHS-based measures of independence and mutual information.

Further work in this vein has shown how RKHS ideas can be used to develop computationally efficient methods for fitting a wide range of other nonparametric and semiparametric models. Consider, for example, the problem of *sufficient dimension reduction* (SDR) for regression (Cook, 1998). This problem can be formulated as the problem of finding a matrix $B_0$ whose column space spans a subspace $\mathcal{S}$ of the covariate space and which satisfies

$$p_{Y|X}(y|x) = p_{Y|B_0^T X}(y|B_0^T x).$$

That is, the regression of $Y$ on $X$ depends only on the projection of $X$ on $\mathcal{S}$. Note that no additional assumptions are made about the regression function. We see that SDR can be viewed in terms of an assertion of *conditional* independence. Fukumizu, Bach and Jordan (2006) have shown how such assertions can be evaluated in terms of covariance operators on RKHSs. In particular, they showed that $B_0$ can be alternatively characterized as

$$B_0 = \arg\min_B \Sigma_{YY|X}^B,$$

where $\Sigma_{YY|X}^B$ is a conditional covariance operator on a RKHS, based on the kernel function $k^B(x, x') = k(B^T x, B^T x')$. This operator can be estimated and minimized (in the sense of the partial order of self-adjoint operators) using Gram matrices. Under weak conditions, this yields a consistent procedure for estimating $B_0$.

## CONCLUSIONS

A statistician who encounters SVMs for the first time might have difficulty understanding the source of the excitement. After all, the SVM is a modest variation on some standard statistical methodology—it involves RKHS expansions of discriminant or regression functions combined with a simple piecewise-linear loss function. Nonetheless, this combination has noteworthy practical consequences. In particular, by paying careful attention to the optimization problem that arises in the SVM and by paying careful attention to the resulting numerical linear algebra, the SVM can be applied to very large classification and regression problems. Moreover, these lessons extend beyond the specific setting of the SVM. As we have emphasized, the key ideas of convex relaxation and reproducing kernels have applications well beyond the SVM. They permit an approach to nonparametric statistics that blends tools from nearby areas of applied mathematics such as optimization theory, functional analysis, numerical linear algebra and combinatorics, undeniably expanding the scope of activities in nonparametric statistics and expanding the scale of problems that can be addressed.